\input amstex
\input amsppt.sty
\magnification=\magstep1
\hsize=30truecc
\baselineskip=16truept
\vsize=22.2truecm
\TagsOnRight
\nologo
\pageno=1
\topmatter
\def\Z{\Bbb Z}
\def\N{\Bbb N}

\def\l{\left}
\def\r{\right}
\def\bg{\bigg}
\def\({\bg(}
\def\[{\bg[}
\def\){\bg)}
\def\]{\bg]}
\def\t{\text}
\def\f{\frac}
\def\mo{\roman{mod}}

\def\bi{\binom}
\def\M#1#2{\thickfracwithdelims[]\thickness0{#1}{#2}_q}
\def\Str#1#2{\left[\matrix#1\\#2\endmatrix\right]}

\def\eq{\equiv}

\def\ls{\leqslant}
\def\gs{\geqslant}
\def\al{\alpha}

\def\da{\delta}

\def\bi{\binom}

\def\Proof{\noindent{\it Proof}}
\def\Remark{\medskip\noindent{\it  Remark}}
\def\Ack{\medskip\noindent {\bf Acknowledgments}}
\hbox{Acta Arith. 124(2006), no.\,1, 41--57.}
\bigskip
\title
On $q$-Euler numbers, $q$-Sali\'e numbers and $q$-Carlitz numbers
\endtitle
\rightheadtext{$q$-Euler numbers, $q$-Sali\'e numbers and $q$-Carlitz numbers}
\author Hao Pan and Zhi-Wei Sun (Nanjing)\endauthor
\leftheadtext{H. Pan and Z. W. Sun}
\abstract Let $(a;q)_n=\prod_{0\ls k<n}(1-aq^k)$ for $n=0,1,2,\ldots$.
Define $q$-Euler numbers $E_{n}(q)$, $q$-Sali\'e numbers $S_{n}(q)$
and $q$-Carlitz numbers $C_n(q)$ as follows:
$$\sum_{n=0}^{\infty}E_{n}(q)\frac{x^{n}}{(q,q)_n}
=\(\sum_{n=0}^{\infty}\frac{q^{n(2n-1)}x^{2n}}{(q;q)_{2n}}\)^{-1},$$
$$\sum_{n=0}^\infty S_{n}(q)\frac{x^{n}}{(q;q)_{n}}
=\sum_{n=0}^\infty \frac{q^{n(n-1)}x^{2n}}{(q;q)_{2n}}
\bigg/\sum_{n=0}^\infty \frac{(-1)^nq^{n(2n-1)}x^{2n}}{(q;q)_{2n}}$$
and
$$\sum_{n=0}^{\infty}C_n(q)\f{x^n}{(q;q)_n}
=\sum_{n=0}^{\infty}\f{q^{n(n-1)}x^{2n+1}}{(q;q)_{2n+1}}
\bigg/\sum_{n=0}^{\infty}\f{(-1)^nq^{n(2n+1)}x^{2n+1}}{(q;q)_{2n+1}}.$$
We show that
$$E_{2n}(q)-E_{2n+2^st}(q)
\equiv [2^s]_{q^t}\ \l(\mo\ (1+q)[2^s]_{q^t}\r)$$
for any nonnegative integers $n,s,t$ with $2\nmid t$,
where $[k]_q=(1-q^k)/(1-q)$; this is a $q$-analogue of
Stern's congruence $E_{2n+2^s}\eq E_{2n}+2^s\ (\mo\ 2^{s+1})$.
We also prove that $(-q;q)_n=\prod_{0<k\ls n}(1+q^k)$ divides $S_{2n}(q)$
and the numerator of $C_{2n}(q)$;
this extends Carlitz's result that $2^n$ divides the Sali\'e number $S_{2n}$ and
the numerator of the Carlitz number $C_{2n}$.
Our result on $q$-Sali\'e numbers implies a conjecture of Guo and Zeng.
\endabstract
\thanks 2000 {\it Mathematics Subject Classifications}:\,Primary 11B65;
Secondary 05A30, 11A07, 11B68.
\newline\indent The second author is responsible for
communications, and supported by the National Science Fund for
Distinguished Young Scholars (no. 10425103) and a Key Program of
NSF in P. R. China.
\endthanks
\endtopmatter
\document

\heading 1. Introduction\endheading

The Euler numbers $E_0,E_1,E_2,\ldots$ are defined by
$$\sum_{n=0}^{\infty}E_n\frac{x^n}{n!}
=\f{2e^x}{e^{2x}+1}=\(\f{e^x+e^{-x}}2\)^{-1}
=\(\sum_{n=0}^{\infty}\f{x^{2n}}{(2n)!}\)^{-1};$$
they are all integers because there holds the recursion
$$\sum^n\Sb k=0\\2\mid k\endSb\bi nkE_{n-k}=\da_{n,0}\quad(n\in\N=\{0,1,2,\ldots\}),$$
where the Kronecker symbol $\da_{n,m}$ is $1$ or $0$
according as $n=m$ or not.
It is easy to see that $E_{2k+1}=0$ for every $k=0,1,2,\ldots$.
In 1871 Stern [St] obtained an interesting arithmetic property of the Euler numbers:
$$E_{2n+2^s}\equiv E_{2n}+2^s\ (\mo\ 2^{s+1})\ \ \t{for any}\ n,s\in\N;
\tag 1.1$$
equivalently we have
$$E_{2m}\equiv E_{2n}\ (\mo\ 2^{s+1})\iff m\equiv n\ (\mo\ 2^s)
\ \ \t{for any}\ m,n,s\in\N.\tag 1.1$'$
$$
Later Frobenius amplified Stern's proof in 1910, and
several different proofs of $(1.1)$ or $(1.1')$
were given by Ernvall [E], Wagstaff [W] and Sun [Su].
Our first goal is to provide a complete $q$-analogue of the Stern congruence.

As usual we let $(a;q)_n=\prod_{0\ls k<n}(1-aq^k)$ for every $n\in\N$,
where an empty product is regarded to have value $1$ and hence $(a;q)_0=1$.
For $n\in\N$ we set
$$[n]_q=\f{1-q^n}{1-q}=\sum_{0\ls k<n}q^k,$$
this is the usual $q$-analogue of $n$. For any $n,k\in\N$,
if $k\ls n$ then we call
$$\M nk=\f{\prod_{0<r\ls n}[r]_q}{(\prod_{0<s\ls k}[s]_q)
(\prod_{0<t\ls n-k}[t]_q)}=\f{(q;q)_n}{(q;q)_k(q;q)_{n-k}}$$
a {\it $q$-binomial coefficient}; if $k>n$ then we let $\M nk=0$.
Obviously we have $\lim_{q\to1}\M nk=\bi nk$. It is easy to see that
$$\M nk=q^k\M{n-1}k+\M{n-1}{k-1}\quad \ \t{for all}\ k,n=1,2,3,\ldots.$$
By this recursion, each $q$-binomial coefficient is a polynomial in $q$ with integer
coefficients.

We define $q$-Euler numbers $E_n(q)\ (n\in\N)$ by
$$\sum_{n=0}^{\infty}E_{n}(q)\frac{x^{n}}{(q;q)_{n}}
=\(\sum_{n=0}^{\infty}\frac{q^{\bi{2n}{2}}x^{2n}}{(q;q)_{2n}}\)^{-1}.\tag1.2$$
Multiplying both sides by $\sum_{n=0}^{\infty}q^{\bi{2n}{2}}x^{2n}/(q;q)_{2n}$, we obtain
the recursion
$$\sum^n\Sb k=0\\2\mid k\endSb \M{n}{k}q^{\bi k2}E_{n-k}(q)=\da_{n,0}\ \ (n\in\N)$$
which implies that $E_n(q)\in\Z[q]$.
Observe that
$$\align&\sum_{n=0}^{\infty}E_n(q)\f{x^n}{\prod_{0<k\ls n}[k]_q}
=\sum_{n=0}^{\infty}E_n(q)\f{((1-q)x)^n}{(q;q)_n}
\\=&\(\sum_{n=0}^{\infty}\frac{q^{\bi{2n}{2}}((1-q)x)^{2n}}{(q;q)_{2n}}\)^{-1}
=\(\sum_{n=0}^{\infty}\frac{q^{\bi{2n}{2}}x^{2n}}{\prod_{0<k\ls 2n}[k]_q}\)^{-1}
\endalign$$
and hence $\lim_{q\to 1}E_n(q)=E_n$.

The usual way to define a $q$-analogue of Euler numbers is as follows:
$$\sum_{n=0}^{\infty}\tilde E_{n}(q)\frac{x^{n}}{(q;q)_{n}}
=\(\sum_{n=0}^{\infty}\frac{x^{2n}}{(q;q)_{2n}}\)^{-1}.$$
(See, e.g., [GZ].) We assert that $\tilde E_n(q)=q^{\bi n2}E_n(1/q)$. In fact,
$$\align&\sum_{n=0}^{\infty}q^{\bi n2}E_n(q^{-1})\f{x^n}{\prod_{0<k\ls n}(1-q^k)}
=\sum_{n=0}^{\infty}E_n(q^{-1})\f{(-q^{-1}x)^n}{\prod_{0<k\ls n}(1-q^{-k})}
\\=&\(\sum_{n=0}^{\infty}\f{q^{-\bi{2n}2}(-q^{-1}x)^{2n}}{\prod_{0<k\ls 2n}(1-q^{-k})}\)^{-1}
=\(\sum_{n=0}^{\infty}\f{x^{2n}}{\prod_{0<k\ls2n}(1-q^k)}\)^{-1}.
\endalign$$

Recently, with the help of cyclotomic polynomials,
Guo and Zeng [GZ] proved that if $m,n,s,t\in\N$, $m-n=2^st$ and $2\nmid t$ then
$$\tilde E_{2m}(q)\equiv q^{m-n}\tilde E_{2n}(q)
\ \(\mo\ \prod_{r=0}^s(1+q^{2^rt})\).$$
This is a partial $q$-analogue of Stern's result.

Using our $q$-analogue of Euler numbers, we are able to give below a complete
$q$-analogue of the classical result of Stern.

\proclaim{Theorem 1.1} Let $n,s,t\in\N$ and $2\nmid t$. Then
$$E_{2n}(q)-E_{2n+2^st}(q)\equiv [2^s]_{q^t}\ \l(\mo\ (1+q)[2^s]_{q^t}\r).\tag1.3$$
\endproclaim

The Sali\'e numbers $S_{n}\ (n\in\N)$ are given by
$$\sum_{n=0}^\infty S_{n}\frac{x^{n}}{n!}
=\frac{\cosh x}{\cos x}=\f{(e^x+e^{-x})/2}{(e^{ix}+e^{-ix})/2}
=\(\sum_{n=0}^\infty \frac{x^{2n}}{(2n)!}\)
\bigg/\sum_{n=0}^\infty (-1)^n\frac{x^{2n}}{(2n)!}.$$
Multiplying both sides by $\sum_{n=0}^\infty (-1)^nx^{2n}/(2n)!$
we get the recursion
$$\sum^n\Sb k=0\\2\mid k\endSb\bi nk (-1)^{k/2}S_{n-k}=\f{1+(-1)^n}2\ \ (n\in\N)$$
which implies that all Sali\'e numbers are integers
and $S_{2k+1}=0$ for all $k\in\N$.

By a sophisticated use of some deep properties of Bernoulli numbers,
in 1965 Carlitz [C2] proved that
$2^n\mid S_{2n}$ for any $n\in\N$ (which was first conjectured by Gandhi [G]).
Recently Guo and Zeng [GZ] defined a $q$-analogue of Sali\'e numbers in the following way:
$$\sum_{n=0}^{\infty}\tilde S_{n}(q)\f{x^{n}}{(q;q)_{n}}
=\sum_{n=0}^{\infty}\f{q^{n^2}x^{2n}}{(q;q)_{2n}}\bigg/
\sum_{n=0}^{\infty}(-1)^n\f{x^{2n}}{(q;q)_{2n}}$$
and hence
$$\sum_{k=0}^n\M{2n}{2k}(-1)^{k}\tilde S_{2n-2k}(q)=q^{n^2}\ \ \t{for any}\ n\in\N.$$
They conjectured that $(-q;q)_n=\prod_{0<k\ls n}(1+q^k)$ divides $\tilde S_{2n}(q)$ (in $\Z[q]$).

We define $q$-Sali\'e numbers by
$$\sum_{n=0}^\infty S_{n}(q)\frac{x^{n}}{(q;q)_{n}}
=\sum_{n=0}^\infty \frac{q^{n(n-1)}x^{2n}}{(q;q)_{2n}}
\bigg/\sum_{n=0}^\infty \frac{(-1)^nq^{\bi{2n}{2}}x^{2n}}{(q;q)_{2n}}.\tag1.4$$
Multiplying both sides by $\sum_{n=0}^\infty (-1)^nq^{\bi{2n}{2}}x^{2n}/(q;q)_{2n}$
one finds the recursion
$$\sum^n_{k=0}\M {2n}{2k} (-1)^{k}q^{\bi {2k}2}S_{2n-2k}(q)=q^{n(n-1)}\ \ \ (n\in\N).\tag1.5$$
In this paper we are able to prove the following $q$-analogue of Carlitz's result
concerning Sali\'e numbers.

\proclaim{Theorem 1.2} Let $n\in\N$. Then
$(-q;q)_n=\prod_{0<k\ls n}(1+q^k)$ divides $S_{2n}(q)$ in the ring $\Z[q]$.
\endproclaim

\proclaim{Corollary 1.1} For any $n\in\N$ we have $(-q;q)_n\mid\tilde S_{2n}(q)$
in the ring $\Z[q]$ as conjectured by Guo and Zeng.
\endproclaim
\Proof. By Theorem 1.2, $S_{2n}(q)=(-q;q)_nP_n(q)$ for some $P_n(q)\in\Z[q]$.
Let $m$ be a natural number not smaller than $\deg P$. Then $q^mP(q^{-1})\in\Z[q]$.
Since
$$q^{\bi{n+1}2}\prod_{0<k\ls n}(1+q^{-k})=\prod_{0<k\ls n}(1+q^k),$$
 $q^{m+\bi{n+1}2}S_{2n}(q^{-1})$ is in $\Z[q]$ and divisible by $(-q;q)_n$.
If the equality
$$\tilde S_{2n}(q)=q^{\bi{2n}2}S_{2n}(q^{-1})$$ holds,
then $q^m\tilde S_{2n}(q)$ is divisible by $(-q;q)_n$ and hence so is $\tilde S_{2n}(q)$
since $q^m$ is relatively prime to $(-q;q)_n$.

Now let us explain why $\tilde S_n(q)=q^{\bi n2}S_n(q^{-1})$ for any $n\in\N$. In fact,
$$\align&\sum_{n=0}^{\infty}q^{\bi n2}S_n(q^{-1})\f{x^n}{\prod_{0<k\ls n}(1-q^k)}
=\sum_{n=0}^{\infty}S_n(q^{-1})\f{(-q^{-1}x)^n}{\prod_{0<k\ls n}(1-q^{-k})}
\\=&\sum_{n=0}^{\infty}\f{q^{-n(n-1)}(-q^{-1}x)^{2n}}{\prod_{0<k\ls 2n}(1-q^{-k})}
\bigg/\sum_{n=0}^{\infty}\f{(-1)^nq^{-\bi{2n}2}(-q^{-1}x)^{2n}}{\prod_{0<k\ls 2n}(1-q^{-k})}
\\=&\sum_{n=0}^{\infty}\f{q^{n^2}x^{2n}}{\prod_{0<k\ls 2n}(1-q^k)}
\bigg/\sum_{n=0}^{\infty}\f{(-1)^nx^{2n}}{\prod_{0<k\ls 2n}(1-q^k)}
=\sum_{n=0}^{\infty}\tilde S_n(q)\f{x^n}{(q;q)_n}.
\endalign$$
This concludes our proof. \qed

In 1956 Carlitz [C1] investigated the coefficients of
$$\f{\sinh x}{\sin x}=\sum_{n=0}^{\infty}C_n\f{x^n}{n!},$$
where
$$\sinh x=\f{e^x-e^{-x}}2=\sum_{n=0}^{\infty}\f{x^{2n+1}}{(2n+1)!}.$$
We call those numbers $C_n\ (n\in\N)$ {\it Carlitz numbers}.
In 1965 Carlitz [C2] proved a conjecture of Gandhi [G] which states that
$2^n$ divides the numerator of $C_{2n}$.

Now we define $q$-Carlitz numbers $C_n(q)\ (n\in\N)$ by
$$\sum_{n=0}^{\infty}C_n(q)\f{x^n}{(q;q)_n}
=\sum_{n=0}^{\infty}\f{q^{n(n-1)}x^{2n+1}}{(q;q)_{2n+1}}
\bigg/\sum_{n=0}^{\infty}\f{(-1)^nq^{n(2n+1)}x^{2n+1}}{(q;q)_{2n+1}}.\tag1.6$$
Multiplying both sides by $\sum_{n=0}^{\infty}(-1)^nq^{n(2n+1)}x^{2n+1}/(q;q)_{2n+1}$
we get the recursion
$$\sum_{k=0}^n\M{2n+1}{2k+1}(-1)^kq^{k(2k+1)}C_{2n-2k}(q)=q^{n(n-1)}\ \ (n\in\N).\tag1.7$$
By (1.7) and induction,
$$[1]_q[3]_q\cdots[2n+1]_qC_{2n}(q)\in\Z[q];$$
in particular, $(2n+1)!!C_{2n}\in\Z$.
If $j,k\in\N$ and $q^j=-1$, then $q^{j(2k+1)}=-1$ and hence $q^{2k+1}\not=1$.
Thus $q^j+1$ is relatively prime to $1-q^{2k+1}$ for any $j,k\in\N$, and hence
$(-q;q)_n=\prod_{0<j\ls n}(1+q^j)$ is relatively prime to
the denominator of $C_{2n}(q)$. This basic property will be used later.

Here is our $q$-analogue of Carlitz's divisibility result concerning Carlitz numbers.

\proclaim{Theorem 1.3} For any $n\in\N$, $(-q;q)_n$ divides the
numerator of $C_{2n}(q)$.
\endproclaim

Note that $E_{2k+1}(q)=S_{2k+1}(q)=C_{2k+1}(q)=0$ for all $k\in\N$ because
$$\sum_{n=0}^{\infty}E_n(q)\f{x^n}{(q;q)_n},
\ \sum_{n=0}^{\infty}S_n(q)\f{x^n}{(q;q)_n},
\ \sum_{n=0}^{\infty}C_n(q)\f{x^n}{(q;q)_n}$$
are even functions.

Our approach to $q$-Euler numbers, $q$-Sali\'e numbers and $q$-Carlitz numbers
is quite different
from that of Guo and Zeng [GZ].
The proofs of Theorems 1.1--1.3 depend on new recursions for
$q$-Euler numbers, $q$-Sali\'e numbers and $q$-Carlitz numbers.
In the next section we will prove Theorem 1.1.
In Section 3 we establish an auxiliary theorem
which essentially says that if $l\in\Z$ and $n\in\N$ then
$$\sum\Sb k\in\Z\\2k+l\gs0\endSb(-1)^kq^{k(k-1)}\M{2n}{2k+l}\eq0\ \ (\mo\ (-q;q)_n).\tag1.8$$
(We can also substitute $2n+1$ for $2n$ in (1.8).)
Section 4 is devoted to the proofs of
Theorems 1.2 and 1.3 on the basis of Section 3.

\heading 2. Proof of Theorem 1.1 \endheading

\proclaim{Lemma 2.1} For any $n\in\N$ we have
$$E_{2n}(q)=1-\sum_{0<k\ls n}(-q;q)_{2k-1}\Str{2n}{2k}_qE_{2(n-k)}(q).\tag 2.1$$
\endproclaim
\Proof. Let us recall the following three known identities
(cf. Theorem 10.2.1 and Corollary 10.2.2 of [AAR]):
$$\sum_{n=0}^{\infty}\frac{q^{\bi{n}{2}}(-x)^n}{(q;q)_n}=(x;q)_\infty$$
where $(x;q)_{\infty}=\prod_{n=0}^{\infty}(1-xq^n)$,
$$\sum_{n=0}^{\infty}\frac{x^n}{(q;q)_n}=\frac{1}{(x;q)_\infty}\ \
\t{and}\ \ \sum_{n=0}^{\infty}\frac{(-1;q)_nx^n}{(q;q)_n}
=\frac{(-x;q)_\infty}{(x;q)_\infty}.$$

Observe that
$$\align\f12\sum_{n=0}^{\infty}E_{n}(q)\frac{x^{n}}{(q;q)_{n}}
=&\(\sum_{n=0}^{\infty}\frac{q^{\bi{n}{2}}x^n}{(q;q)_n}
+\sum_{n=0}^{\infty}\frac{q^{\bi{n}{2}}(-x)^n}{(q;q)_n}\)^{-1}
\\=&\frac{1}{(x;q)_\infty+(-x;q)_\infty}
\endalign$$
and hence
$$\align&\f12\(\sum_{n=0}^{\infty}E_n(q)\f{x^n}{(q;q)_n}\)
\(1+\sum_{n=0}^{\infty}\f{(-1;q)_nx^n}{(q;q)_n}\)
\\=&\f1{(x;q)_{\infty}+(-x;q)_{\infty}}\(1+\f{(-x;q)_{\infty}}{(x;q)_{\infty}}\)
=\f1{(x;q)_{\infty}}=\sum_{n=0}^{\infty}\f{x^n}{(q;q)_n}.
\endalign$$
Comparing the coefficients of $x^n$ we obtain that
$$\f12E_n(q)+\f12\sum_{k=0}^n(-1;q)_k\M nk E_{n-k}(q)=1,$$
i.e.,
$$E_n(q)=1-\sum_{0<k\ls n}(-q;q)_{k-1}\M nk E_{n-k}(q).$$
Substituting $2n$ for $n$ in the last equality
and recalling that $E_{2j+1}(q)=0$ for $j\in\N$, we immediately
obtain the desired (2.1). \qed

\proclaim{Corollary 2.1} For any $n\in\N$ we have
$$E_{2n}(q)\equiv1\ (\mo\ 1+q).\tag2.2$$
\endproclaim
\Proof. This follows from $(2.1)$ because $1+q$ divides
$(-q;q)_{m}$ for all $m=1,2,3,\ldots$. \qed
\medskip

The following trick is simple but useful.
$$\prod_{k=0}^n(1+q^{2^k})=[2^{n+1}]_q\ \ \ \t{for any}\ n\in\N.\tag2.3$$
In fact,
$$\align (1-q)\prod_{k=0}^n(1+q^{2^k})=&(1-q^2)\prod_{0<k\ls n}(1+q^{2^k})
\\=&\cdots=(1-q^{2^n})(1+q^{2^n})=1-q^{2^{n+1}}.
\endalign$$

\proclaim{Lemma 2.2} Let $m,n,s,t$ be positive
integers with $2m\gs n$ and $2\nmid t$. Then
$(-q;q)_m\M{2^st}{n}$ is divisible by
$(1+q)^{\lfloor(m-1)/2\rfloor}[2^s]_{q^t}$, where we use $\lfloor \al\rfloor$
to denote the greatest integer not exceeding a real number $\al$.
\endproclaim
\Proof. Write $n=2^kl$ with $k,l\in\N$ and $2\nmid l$. Then
$$[n]_q=\frac{1-q^n}{1-q}=\frac{1-q^{2^kl}}{1-q^l}\cdot\frac{1-q^l}{1-q}
=[2^k]_{q^l}[l]_q.$$
Obviously $[2^k]_{q^l}=\prod_{0\ls j<k}(1+q^{2^jl})$ divides
$(-q;q)_{m}=\prod_{j=1}^{m}(1+q^j)$ since $m\gs n/2=2^{k-1}l$.
Thus $[2^s]_{q^t}=[2^st]_q/[t]_q$ divides
$$[l]_q(-q;q)_m\M{2^st}{n}=\frac{(-q;q)_m}{[2^k]_{q^l}}
[2^st]_{q}\M{2^st-1}{n-1}.$$
Note that $[2^s]_{q^t}=\prod_{r=0}^{s-1}(1+q^{2^rt})$
is relatively prime to $[l]_q=(1-q^l)/(1-q)$ since $l\eq1\ (\mo\ 2)$.
Therefore $[2^s]_{q^t}$ divides $(-q;q)_m\M{2^st}{n}$.

Clearly $(1+q)^{\lfloor(m+1)/2\rfloor}$ divides
$$\prod_{j=1}^{\lfloor(m+1)/2\rfloor}(1+q^{2j-1})
\times\prod_{j=1}^{\lfloor m/2\rfloor}(1+q^{2j})=(-q;q)_m.$$
Since
$$[2^s]_{q^t}=\f{1-q^{2t}}{1-q^t}\cdot\f{1-q^{2^st}}{1-q^{2t}}
=(1+q)\sum_{j=0}^{t-1}(-q)^j\sum_{r=0}^{2^{s-1}-1}q^{2rt}$$
and $\sum_{0\ls j<t}(-q)^j\sum_{0\ls r<2^{s-1}}q^{2rt}$ takes value $2^{s-1}t\not=0$ at $q=-1$,
the polynomial $[2^s]_{q^t}$ is divisible by $1+q$ but not by $(1+q)^2$.
Therefore
$(1+q)^{\lfloor (m-1)/2\rfloor}[2^s]_{q^t}$ divides $(-q;q)_m\M{2^st}n$
by the above.  \qed

\medskip

\noindent{\it Proof of Theorem 1.1}. The case $s=0$ is easy. In fact,
$$E_{2n}(q)-E_{2n+2^0t}(q)=E_{2n}(q)\eq1=[2^0]_{q^t}\ (\mo\ (1+q)[2^0]_{q^t})$$
by Corollary 2.1.

Below we handle the case $s>0$ and use induction on $n$.
Assume that
$$E_{2m}(q)-E_{2m+2^st}(q)\equiv [2^s]_{q^t}\ (\mo\ (1+q)[2^s]_{q^t})\
\t{whenever}\ 0\ls m<n.\tag $*$ $$
(This holds trivially in the case $n=0$.)
In view of Lemma 2.1, we have
$$\aligned
&E_{2n}(q)-E_{2n+2^st}(q)\\
=&\sum_{k=1}^{n+2^{s-1}t}(-q;q)_{2k-1}\(\M{2n+2^st}{2k}
E_{2n+2^st-2k}(q)-\M{2n}{2k}E_{2n-2k}(q)\),
\endaligned$$
where we set $E_{l}(q)=0$ for $l<0$.

Let $0<k\ls n+2^{s-1}t$.
Applying a
$q$-analogue of the Chu-Vandermonde identity (cf. [AAR, Exercise 10.4(b)]), we find that
$$\aligned
&\M{2n+2^st}{2k}E_{2n+2^st-2k}(q)-\M{2n}{2k}E_{2n-2k}(q)
\\=&E_{2n+2^st-2k}(q)\sum_{j=0}^{2k}q^{(2n-j)(2k-j)}\M{2n}{j}\M{2^st}{2k-j}
-\M{2n}{2k}E_{2n-2k}(q)
\\=&E_{2n+2^st-2k}(q)\sum_{j=0}^{2k-1}q^{(2n-j)(2k-j)}\M{2n}{j}\M{2^st}{2k-j}
\\&+\M{2n}{2k}(E_{2n+2^st-2k}(q)-E_{2n-2k}(q)).
\endaligned$$
In view of the hypothesis $(*)$,
$$(-q;q)_{2k-1}\M{2n}{2k}(E_{2n+2^st-2k}(q)-E_{2(n-k)}(q))\equiv 0\ (\mo\ (1+q)[2^s]_{q^t}).$$
By Lemma 2.2, if $0\ls j<2k$ then
$(-q;q)_{2k-1}\M{2^st}{2k-j}$ is divisible by $(1+q)^{k-1}[2^s]_{q^t}$.
Therefore, if $k>1$ then
$(1+q)[2^s]_{q^t}$ divides
$$(-q;q)_{2k-1}\(\M{2n+2^st}{2k}E_{2n+2^st-2k}(q)-\M{2n}{2k}E_{2n-2k}(q)\)$$
by the above. In the case $k=1$,
$$(-q;q)_{2k-1}\M{2^st}{2k-1}=(1+q)[2^st]_q=(1+q)[2^s]_{q^t}[t]_q$$
and hence
$$\align&(-q;q)_1\(\M{2n+2^st}{2}E_{2n+2^st-2}(q)-\M{2n}2E_{2n-2}(q)\)
\\\eq&(1+q)E_{2n+2^st-2}(q)q^{(2n-0)(2-0)}\M{2n}0\M{2^st}2\ (\mo\ (1+q)[2^s]_{q^t})
\\\eq&E_{2n+2^st-2}(q)q^{4n}\f{1+q}{[2]_q}[2^st]_q[2^st-1]_q\  (\mo\ (1+q)[2^s]_{q^t})
\\\eq&E_{2n+2^st-2}(q)q^{4n}[2^s]_{q^t}[t]_q(1+q[2^st-2]_q)
\eq[2^s]_{q^t}\ (\mo\ (1+q)[2^s]_{q^t});
\endalign$$
in the last step we have noted that $q^{4n}-1,[t]_q-1, [2^st-2]_q$ are
divisible by $1+q$,
and $E_{2n+2^st-2}(q)\eq 1\ (\mo\ 1+q)$ by Corollary 2.1.

Combining the above we obtain that
$$E_{2n}(q)-E_{2n+2^st}(q)
\eq\sum_{k=1}^{n+2^{s-1}t}\da_{k,1}[2^s]_{q^t}=[2^s]_{q^t}\ (\mo\ (1+q)[2^s]_{q^t}).$$
This concludes the induction.

 The proof of Theorem 1.1 is now complete. \qed

 \Remark\ 2.1. With a bit more efforts
 we can prove the following more general result:
For $k=1,2,3,\ldots$ let
$$\sum_{n=0}^{\infty}E_n^{(k)}(q)\f{x^n}{(q;q)_n}
=\(\sum_{n=0}^{\infty}q^{\bi{kn}2}\f{x^{kn}}{(q;q)_{kn}}\)^{-1}.$$
Given positive integers $k,s,t$ with $2\nmid t$, we have
$$E^{(2k')}_{2k'n}(q)-E^{(2k')}_{2k'(n+2^{s-1}t)}(q)
\eq(2k'-1)[2^s]_{q^{k't}}
\ \l(\mo\ (1+q^{k'})[2^s]_{q^{k't}}\r)$$
for all $n\in\N$, where $k'=2^{k-1}$.
This is a $q$-analogue of Conjecture 5.5 in [GZ].

 \heading {3. An Auxiliary Theorem}\endheading

\proclaim{Theorem 3.1} For all $m,n\in\N$, both
$$S^m_n:=\sum_{k=0}^n(-1)^kq^{k(k-1)+2m(n-k)}\M{2n}{2k}\tag3.1$$
and
$$T^m_n:=\sum_{0\ls k<n}(-1)^kq^{k(k-1)+2m(n-1-k)}\M{2n}{2k+1}\tag3.2$$
are divisible by $(-q;q)_n=\prod_{0<k\ls n}(1+q^k)$ in the ring $\Z[q]$.
Also, for any $m,n\in\N$ and $\da\in\{0,1\}$ we have the congruence
$$\sum_{k=0}^n(-1)^kq^{k(k+2m-1)}\M{2n}{2k+\da}\eq0\ \ (\mo\ (-q;q)_n).\tag3.3$$
\endproclaim
\Proof. (i) We use induction on $n$ to prove the first part.

For any $m\in\N$, clearly both $S^m_0=1$ and $T^m_0=0$
are divisible by $(-q;q)_0=1$, also both
$S^m_1=q^{2m}-1$ and $T^m_1=[2]_q=1+q$ are multiples of $(-q;q)_1=1+q$.

Now let $n>1$ be an integer and assume that
$(-q;q)_{n-1}$ divides both $S^m_{n-1}$ and $T^m_{n-1}$ for all $m\in\N$.

For each $m\in\Z$ we have
$$\align S^m_n=&\sum_{l=0}^n(-1)^{n-l}q^{(n-l)(n-l-1)+2ml}\M{2n}{2(n-l)}
\\=&(-1)^nq^{n(n-1)}\sum_{l=0}^n(-1)^lq^{l(l+1)-2ln+2lm}\M{2n}{2l}
\\=&(-1)^nq^{n(n-1)-2n(n-1-m)}S^{n-1-m}_n=(-1)^nq^{n(2m-n+1)}S^{n-1-m}_n.
\endalign$$
In particular,
$$S_n^n=(-1)^nq^{n(n+1)}S^{-1}_n\ \ \t{and}\ \ S_n^{n-1}=(-1)^nq^{n(n-1)}S^0_n.$$
Similarly, for every $m\in\Z$ we have
$$\align T^m_n=&\sum_{l=0}^{n-1}(-1)^{n-1-l}q^{(n-1-l)(n-l-2)+2ml}\M{2n}{2(n-1-l)+1}
\\=&(-1)^{n-1}q^{(n-1)(n-2)}\sum_{l=0}^{n-1}(-1)^lq^{l(l+1)-2l(n-1)+2lm}\M{2n}{2l+1}
\\=&(-1)^{n-1}q^{(n-1)(2m-n+2)}T^{n-2-m}_n.
\endalign$$
In particular,
$$T^{n-1}_n=(-1)^{n-1}q^{n(n-1)}T_n^{-1}
\ \ \t{and}\ \ T_n^{n-2}=(-1)^{n-1}q^{(n-1)(n-2)}T^0_n.$$

For any $m\in\N$, clearly
$$\align S^{m+1}_n-S^m_n=&\sum_{k=0}^n(-1)^kq^{k(k-1)+2m(n-k)}(q^{2(n-k)}-1)\M{2n}{2k}
\\=&\sum_{k=0}^n(-1)^kq^{k(k-1)+2m(n-k)}(q^{2n}-1)\M{2n-1}{2k}
\\=&(q^{2n}-1)\sum_{k=0}^{n-1}(-1)^kq^{k(k-1)+2m(n-k)}q^{2k}\M{2n-2}{2k}
\\&+(q^{2n}-1)\sum_{k=1}^{n-1}(-1)^kq^{k(k-1)+2m(n-k)}\M{2n-2}{2k-1}
\\=&(q^{2n}-1)q^{2(m+n-1)}S^{m-1}_{n-1}-(q^{2n}-1)q^{2(m+n-2)}T^{m-1}_{n-1}
\\=&(q^{2n}-1)q^{2(m+n-2)}(q^2S^{m-1}_{n-1}-T^{m-1}_{n-1})
\endalign$$
and
$$\align qT^{m+1}_n-T^m_n=&\sum_{k=0}^{n-1}(-1)^kq^{k(k-1)+2m(n-1-k)}
(q^{2(n-1-k)+1}-1)\M{2n}{2k+1}
\\=&\sum_{k=0}^{n-1}(-1)^kq^{k(k-1)+2m(n-1-k)}(q^{2n}-1)\M{2n-1}{2k+1}
\\=&(q^{2n}-1)\sum_{k=0}^{n-2}(-1)^kq^{k(k-1)+2m(n-1-k)}q^{2k+1}\M{2n-2}{2k+1}
\\&+(q^{2n}-1)\sum_{k=0}^{n-1}(-1)^kq^{k(k-1)+2m(n-1-k)}\M{2n-2}{2k}
\\=&(q^{2n}-1)q^{2m+2n-3}T^{m-1}_{n-1}+(q^{2n}-1)S^{m}_{n-1},
\endalign$$
therefore by the induction hypothesis
we have
$$S^{m+1}_n\eq S^m_n\ (\mo\ (-q;q)_n) \ \ \t{and}\ \ qT^{m+1}_n\eq T^m_n\ (\mo\ (-q;q)_n).$$
(Note that both $q^{n(n-1)}S_{n-1}^{-1}=(-1)^{n-1}S^{n-1}_{n-1}$ and
$q^{(n-1)(n-2)}T_{n-1}^{-1}=(-1)^{n}T_{n-1}^{n-2}$ are divisible by $(-q;q)_{n-1}$
by the induction hypothesis.)
Thus, if $(-q;q)_n$ divides both $S^0_n$ and $T^0_n$ then it divides
both $S^m_n$ and $T^m_n$ for every $m=0,1,2,\ldots$.

Observe that
$$\align &S_n^0=\sum_{k=0}^n(-1)^kq^{k(k-1)}\M{2n}{2n-2k}
\\=&\sum_{k=1}^n(-1)^kq^{k(k-1)+2n-2k}\M{2n-1}{2n-2k}
+\sum_{k=0}^{n-1}(-1)^kq^{k(k-1)}\M{2n-1}{2n-2k-1}
\\=&\sum_{k=1}^n(-1)^kq^{k(k-1)}q^{2(2n-2k)}\M{2n-2}{2n-2k}
\\&+\sum_{k=1}^{n-1}(-1)^kq^{k(k-1)}(q^{2n-2k}+q^{2n-2k-1})\M{2n-2}{2n-2k-1}
\\&+\sum_{k=0}^{n-1}(-1)^kq^{k(k-1)}\M{2n-2}{2n-2k-2}
\\=&-q^{2n-2}S_{n-1}^1-q^{2n-3}(1+q)T_{n-1}^0+S_{n-1}^0
\endalign$$
and hence $(-q;q)_{n-1}$ divides $S_n^0$ by the induction hypothesis.
Similarly, $(-q;q)_{n-1}$ divides $T_n^0=-q^{2n-2}T_{n-1}^1+(1+q)S^1_{n-1}+T^0_{n-1}$.

Since
$$(-1)^nq^{n(n-1)}S_n^0=S_n^{n-1}\eq S_n^0\ (\mo\ (-q;q)_n)$$
and
$$1-(-1)^nq^{n(n-1)}\eq1-(-1)^n(-1)^{n-1}=2\ (\mo\ 1+q^n),$$
we must have
$S_n^0/(-q;q)_{n-1}\eq0\ (\mo\ 1+q^n)$ and hence $(-q;q)_n\mid S^0_n$.
Similarly, as
$$q^{n-2}(-1)^{n-1}q^{(n-1)(n-2)}T^0_n=q^{n-2}T^{n-2}_n\eq T^0_n\ (\mo\ (-q;q)_n)$$
and $1-(-1)^{n-1}q^{n(n-2)}\eq2\ (\mo\ 1+q^n)$, we have
$T_n^0/(-q;q)_{n-1}\eq0\ (\mo\ 1+q^n)$ and hence $(-q;q)_n\mid T^0_n$.
This concludes our induction step and proves the first part.

(ii) Now fix $m,n\in\N$ and $\da\in\{0,1\}$.
We can verify (3.3) directly if $n<2$.

Below we assume $n\gs 2$.  By a previous argument,
$$ (-1)^nS_n^{m+n-1}=q^{n(2m+n-1)}S_n^{-m}=q^{n(n-1)}
\sum_{k=0}^n(-1)^kq^{k(k+2m-1)}\M{2n}{2k}$$
and
$$\align(-1)^{n-1}T_n^{m+n-2}
=&q^{(n-1)(2m+n-2)}T_n^{-m}
\\=&q^{(n-1)(n-2)}\sum_{k=0}^{n-1}
(-1)^kq^{k(k+2m-1)}\M{2n}{2k+1}.
\endalign$$
Thus, applying the first part we immediately get (3.3).

The proof of Theorem 3.1 is now complete. \qed

\Remark\ 3.1. Theorem 3.1 is somewhat difficult and sophisticated,
however it is easy to evaluate the sums
$$\sum_{k=0}^n(-1)^k\bi{2n}{2k}=\sum_{k=0}^{2n}\bi{2n}k\f{i^k+(-i)^k}2$$
and
$$\sum_{0\ls k<n}(-1)^k\bi{2n}{2k+1}=\sum_{k=0}^{2n}\bi{2n}k\f{i^k-(-i)^k}{2i}.$$

\medskip

Now let us explain why (1.8) holds for any $l\in\Z$ and $n\in\N$.
Write $l=2m+\da$ with $m\in\Z$ and $\da\in\{0,1\}$. Then
$$\align&\sum\Sb k\in\Z\\2k+l\gs0\endSb(-1)^kq^{k(k-1)}\M{2n}{2k+l}
\\=&\sum_{k+m\in\N}(-1)^kq^{k(k-1)}\M{2n}{2(k+m)+\da}
\\=&\sum_{k\in\N}(-1)^{k-m}q^{(k-m)(k-m-1)}\M{2n}{2k+\da}
\\=&(-1)^m\sum_{k=0}^{n-\da}(-1)^kq^{k(k-1)-2km+m(m+1)}\M{2n}{2k+\da}.
\endalign$$
So (1.8) follows from Theorem 3.1. Note also that
$$\align&\sum\Sb k\in\Z\\2k+l\gs0\endSb(-1)^kq^{k(k-1)}\M{2n+1}{2k+l}
\\&-\sum\Sb k\in\Z\\2k+l-1\gs0\endSb(-1)^kq^{k(k-1)}\M{2n}{2k+l-1}
\\=&\sum\Sb k\in\Z\\2k+l\gs0\endSb(-1)^kq^{k(k-1)+2k+l}\M{2n}{2k+l}
\\=&q^l\sum\Sb k\in\Z\\2k+l-2\gs0\endSb(-1)^{k-1}q^{k(k-1)}\M{2n}{2k+l-2}
\endalign$$
and thus
$$\sum\Sb k\in\Z\\2k+l\gs0\endSb(-1)^kq^{k(k-1)}\M{2n+1}{2k+l}
\eq0\ \ (\mo\ (-q;q)_n).\tag3.4$$

\heading 4. Proofs of Theorems 1.2 and 1.3\endheading
\proclaim{Lemma 4.1} We have
$$1+\sum_{n=1}^\infty (-q;q)_{2n-1}\frac{(-1)^nx^{2n}}{(q;q)_{2n}}
=\sum_{k=0}^\infty q^{\bi{2k}{2}}\frac{(-1)^kx^{2k}}{(q;q)_{2k}}
\sum_{l=0}^\infty \frac{(-1)^lx^{2l}}{(q;q)_{2l}}\tag4.1$$
and
$$\sum_{n=0}^\infty (-q;q)_{2n}\frac{(-1)^nx^{2n+1}}{(q;q)_{2n+1}}
=\sum_{k=0}^\infty q^{\bi{2k+1}{2}}\frac{(-1)^kx^{2k+1}}{(q;q)_{2k+1}}
\sum_{l=0}^\infty \frac{(-1)^lx^{2l}}{(q;q)_{2l}}.\tag4.2$$
\endproclaim
\Proof. Let $\da\in\{0,1\}$. Then
$$\align&\sum_{k=0}^\infty \frac{(-1)^kq^{\bi{2k+\da}{2}}x^{2k+\da}}{(q;q)_{2k+\da}}
\sum_{l=0}^\infty \frac{(-1)^lx^{2l}}{(q;q)_{2l}}
\\=&\sum_{n=0}^\infty\frac{(-1)^nx^{2n+\da}}{(q;q)_{2n+\da}}
\sum_{k=0}^nq^{\bi{2k+\da}{2}}\M{2n+\da}{2k+\da}.
\endalign$$
By the $q$-binomial theorem (cf. [AAR, Corollary 10.2.2(c)]),
$$(x;q)_m=\sum_{k=0}^m\M mk(-1)^kq^{\bi k2}x^k\ \ \t{for any}\ m\in\N.$$
Thus
$$\align2\sum_{k=0}^nq^{\bi{2k+\da}{2}}\M{2n+\da}{2k+\da}
=&\sum_{l=0}^{2n+\da}q^{\bi{l}{2}}\M{2n+\da}{l}
+\sum_{l=0}^{2n+\da}(-1)^{\da+l}q^{\bi{l}{2}}\M{2n+\da}{l}
\\=&(-1;q)_{2n+\da}+(-1)^{\da}(1;q)_{2n+\da}
\endalign$$
and hence
$$\aligned
&\sum_{k=0}^\infty \frac{(-1)^kq^{\bi{2k+\da}{2}}x^{2k+\da}}{(q;q)_{2k+\da}}
\sum_{l=0}^\infty \frac{(-1)^lx^{2l}}{(q;q)_{2l}}
\\=&\sum_{n=0}^\infty\frac{(-1)^nx^{2n+\da}}{(q;q)_{2n+\da}}
\cdot\frac{(-1;q)_{2n+\da}+(-1)^{\da}(1;q)_{2n+\da}}{2}
\\=&\cases1+\sum_{n=1}^\infty(-q;q)_{2n-1}\frac{(-1)^nx^{2n}}{(q;q)_{2n}}&\t{if}\ \da=0,
\\\sum_{n=0}^\infty(-q;q)_{2n}\frac{(-1)^nx^{2n+1}}{(q;q)_{2n+1}}&\t{if}\ \da=1.\endcases
\endaligned$$
We are done. \qed

\Remark\ 4.1. (4.1) and (4.2) are $q$-analogues of the trigonometric identities
$$\f{1+\cos(2x)}2=\cos^2x\ \ \ \t{and}\ \ \ \f{\sin(2x)}2=\sin x\cos x$$
respectively.

\proclaim{Lemma 4.2} Let $n\gs k\gs1$ be integers.
Then both $(-q;q)_{k}\M{2n}{2k}$ and $(-q;q)_{k}\M{2n+1}{2k+1}$ are divisible by
$$(-q^{n-k+1};q)_{k}=\prod_{j=1}^{k}(1+q^{n-j+1}).$$
\endproclaim
\Proof. Observe that
$$\align
\M{2n}{2k}=&\prod_{j=1}^{2k}\f{1-q^{2n-j+1}}{1-q^j}
=\prod_{j=1}^k\f{(1-q^{2n-2j+1})(1-q^{2n-(2j-1)+1})}{(1-q^{2j})(1-q^{2j-1})}
\\=&\prod_{j=1}^k\f{(1-q^{n-j+1})(1+q^{n-j+1})(1-q^{2n-2j+1})}
{(1-q^j)(1+q^j)(1-q^{2j-1})}
\\=&\M nk\f{\prod_{j=1}^k(1+q^{n-j+1})}{(-q;q)_k}\prod_{j=1}^k\f{1-q^{2n-2j+1}}{1-q^{2j-1}}
\endalign$$
and hence
$$(-q^{n-k+1};q)_k\ \bigg|\ (-q;q)_k\M{2n}{2k}\prod_{j=1}^k(1-q^{2j-1}).$$
Recall that
$(-q^{n-k+1};q)_k=\prod_{n-k<l\ls n}(1+q^l)$ is relatively prime to
$\prod_{j=1}^k(1-q^{2j-1})$. Therefore $(-q^{n-k+1};q)_k\mid
(-q;q)_k\M{2n}{2k}$.

Since $[2k+1]_q$ is also relatively prime to $(-q^{n-k+1};q)_k$, we have
$$(-q;q)_k\M{2n+1}{2k+1}=(-q;q)_k\f{[2n+1]_q}{[2k+1]_q}\M{2n}{2k}
\eq0\ (\mo\ (-q^{n-k+1};q)_k).$$
This concludes the proof. \qed

\Remark\ 4.2. Lemma 4.2 yields a trivial result as $q\to 1$.

\medskip
\noindent{\it Proof of Theorem 1.2}. Clearly
$$\align &f(x):=\(\sum_{n=0}^\infty S_{2n}(q)\frac{x^{2n}}{(q;q)_{2n}}\)
\(1+\sum_{n=1}^\infty (-q;q)_{2n-1}\frac{(-1)^nx^{2n}}{(q;q)_{2n}}\)
\\=&\sum_{n=0}^\infty S_{2n}(q)\frac{x^{2n}}{(q;q)_{2n}}
+\sum_{n=1}^\infty \frac{x^{2n}}{(q;q)_{2n}}\sum_{k=1}^n(-1)^k(-q;q)_{2k-1}
\M{2n}{2k}S_{2n-2k}(q).
\endalign$$
On the other hand, by (4.1) we have
$$\align f(x)=&\sum_{k=0}^\infty \frac{q^{k(k-1)}x^{2k}}{(q;q)_{2k}}
\sum_{l=0}^\infty \frac{(-1)^lx^{2l}}{(q;q)_{2l}}
\\=&\sum_{n=0}^\infty \frac{(-1)^nx^{2n}}{(q;q)_{2n}}\sum_{k=0}^n(-1)^kq^{k(k-1)}\M{2n}{2k}.
\endalign$$
Therefore
$$\align&S_{2n}(q)+\sum_{0<k\ls n}(-1)^k(-q;q)_{2k-1}\M{2n}{2k}S_{2n-2k}(q)
\\=&(-1)^n\sum_{k=0}^n(-1)^kq^{k(k-1)}\M{2n}{2k}\eq0\ \ (\mo\ (-q;q)_n)
\endalign$$
with the help of (1.8) or Theorem 3.1.
If $(-q;q)_l\mid S_{2l}(q)$ for all $0\ls l<n$, then
$$S_{2n}(q)\eq-\sum_{0<k\ls n}(-1)^k(-q;q)_{2k-1}\M{2n}{2k}S_{2n-2k}(q)\eq0\ (\mo\ (-q;q)_n)$$
since $\prod_{0<j\ls n-k}(1+q^j)$ divides $S_{2n-2k}(q)$ and
$\prod_{n-k<j\ls n}(1+q^j)$ divides $(-q;q)_{2k-1}\M{2n}{2k}$ by Lemma 4.2.
Thus we have the desired result by induction. \qed

\Remark\ 4.3. As $q\to1$ our new recursion for $q$-Sali\'e numbers yields
a useful recursion for Sali\'e numbers:
$$S_{2n}+\sum_{0<k\ls n}(-1)^k2^{2k-1}\bi{2n}{2k}S_{2n-2k}=(-1)^n\sum_{k=0}^n(-1)^k\bi{2n}{2k},$$
from which the Carlitz result $2^n\mid S_{2n}$ follows by induction.

\medskip
\noindent{\it Proof of Theorem 1.3}. It is apparent that
$$\align g(x):=&\(\sum_{n=0}^\infty C_{2n}(q)\frac{x^{2n}}{(q;q)_{2n}}\)
\(\sum_{n=0}^\infty (-q;q)_{2n}\frac{(-1)^nx^{2n+1}}{(q;q)_{2n+1}}\)
\\=&\sum_{n=0}^\infty \frac{x^{2n+1}}{(q;q)_{2n+1}}\sum_{k=0}^n(-1)^k(-q;q)_{2k}
\M{2n+1}{2k+1}C_{2n-2k}(q).
\endalign$$
On the other hand, (4.2) implies that
$$\align g(x)=&\sum_{k=0}^\infty \frac{q^{k(k-1)}x^{2k+1}}{(q;q)_{2k+1}}
\sum_{l=0}^\infty \frac{(-1)^lx^{2l}}{(q;q)_{2l}}
\\=&\sum_{n=0}^\infty \frac{x^{2n+1}}{(q;q)_{2n+1}}
\sum_{k=0}^n(-1)^{n-k}q^{k(k-1)}\M{2n+1}{2k+1}.
\endalign$$
Therefore we have the recurrence relation
$$\sum_{k=0}^n(-1)^k(-q;q)_{2k}\M{2n+1}{2k+1}C_{2n-2k}(q)
=\sum_{k=0}^n(-1)^{n-k}q^{k(k-1)}\M{2n+1}{2k+1}.$$
The right-hand side of the last equality is a multiple of $(-q;q)_n$ by (3.4).
So we have
$$\sum_{k=0}^n(-1)^k(-q;q)_{2k}\M{2n+1}{2k+1}C_{2n-2k}(q)\eq0\ \ (\mo\ (-q;q)_n).$$

Assume that $(-q;q)_l$ divides the numerator of $C_{2l}(q)$ for each $0\ls l<n$. Then
$(-q;q)_n$ divides the numerator of $(-q;q)_{2k}\M{2n+1}{2k+1}C_{2n-2k}(q)$
for each $0<k\ls n$,
because $\prod_{0<j\ls n-k}(1+q^j)$ divides the numerator of $C_{2n-2k}(q)$ and
$\prod_{n-k<j\ls n}(1+q^j)$ divides $(-q;q)_{2k}\M{2n+1}{2k+1}$ by Lemma 4.2.
Thus $(-q;q)_n$ must also divide the numerator of
$\M{2n+1}1C_{2n}(q)=[2n+1]_qC_{2n}(q)$.
Recall that $[2n+1]_q$ is relatively prime to $(-q;q)_n$. So the numerator of $C_{2n}(q)$
is divisible by $(-q;q)_n$.

In view of the above, the desired result follows by induction on $n$.
\qed

\Remark\ 4.4. As $q\to1$ our new recursion for $q$-Carlitz numbers yields
the following recurrence relation for Carlitz numbers:
$$\sum_{k=0}^n(-1)^k2^{2k}\bi{2n+1}{2k+1}C_{2n-2k}=(-1)^n\sum_{k=0}^n(-1)^k\bi{2n+1}{2k+1}.$$
From this one can easily deduce the Carlitz congruence $C_{2n}\eq0\ (\mo\ 2^n)$.

\Ack. The second author is indebted to Prof. Jiang Zeng
at University of Lyon-I for showing
Carlitz's paper [C2] and the preprint [GZ]
during Sun's visit to the Institute of Camille Jordan.
This paper was finished during Sun's visit to the University of California
at Irvine; he would like to thank Prof. Daqing Wan for the invitation.

\bigskip

\leftline{Department of Mathematics}
\leftline{Nanjing University}
\leftline{Nanjing 210093}
\leftline{People's Republic of China}
\leftline {E-mail: (Hao Pan) {\tt haopan79\@yahoo.com.cn}}
\leftline {\qquad (Zhi-Wei Sun) {\tt zwsun\@nju.edu.cn}}


\medskip

\widestnumber\key{AAR}

\Refs

\ref\key AAR\by G. E. Andrews, R. Askey and R. Roy\book Special Functions
\publ Cambridge University Press, Cambridge, 1999\endref

\ref\key C1\by L. Carlitz\paper The coefficients of $\sinh x/\sin x$
\jour Math. Mag.\vol 29\yr1956\pages 193-197\endref

\ref\key C2\by L. Carlitz\paper The coefficients of $\cosh x/\cos x$
\jour Monatsh. Math.\vol69\yr1965\pages 129-135\endref

\ref\key E\by R. Ernvall\paper Generalized Bernoulli numbers, generalized irregular primes,
 and class number\jour Ann. Univ. Turku. Ser. A, I(178), 1979, 72 pp\endref

\ref\key G\by J. M. Gandhi\paper The coefficients of $\cosh x/\cos x$
and a note on Carlitz's coefficients of $\sinh x/\sin x$
\jour Math. Mag.\vol 31\yr1958\pages 185-191\endref

\ref\key GZ\by Victor J. W. Guo and J. Zeng\paper Some arithmetic
properties of the $q$-Euler numbers and q-Sali\'e numbers \jour
European J. Combin.\vol 27\yr 2006\pages 884--895\endref

\ref\key St\by M. A. Stern\paper Zur Theorie der Eulerschen Zahlen
\jour J. Reine Angew. Math.\vol79\yr1875\pages 67-98\endref

\ref\key Su\by Z. W. Sun\paper On Euler numbers modulo powers of
two\jour J. Number Theory\vol 115\yr 2005\pages 371--380\endref

\ref\key W\by S. S. Wagstaff, Jr.\paper Prime divisors of the Bernoulli and Euler numbers
\jour in: Number Theory for the Millennium, III (Urbana, IL, 2000), 357--374,
A K Peters, Natick, MA, 2002\endref

\endRefs
\enddocument